# Enhanced Ramping Product Design to Improve Power System Reliability

Hyeongon Park, *Member, IEEE*, Bing Huang, *Student Member, IEEE*, and Ross Baldick, *Fellow, IEEE*

*Abstract*— The roll-out of a flexible ramping product provides Independent System Operators (ISOs) with the ability to address ramping capacity shortages. ISOs procure flexible ramping capability by committing more generating units or reserving a certain amount of head room of committed units. In this paper, we raise the concern of the possibility that the procured ramping capability cannot be deployed in real-time operations. As a solution to the non-delivery issue, we provide a new ramping product designed to improve reliability and reduce the expected operating cost. The trajectory of start-up and shutdown processes is also considered in determining the ramping capability. A new optimization problem is formulated using mixed integer linear programming to be readily applied to practical power system operation. The performance of this proposed method is verified through simulations using a small-scale system and IEEE 118-bus systems. The simulation results demonstrate that the proposed formulation can realize improved generation scheduling alleviating capacity shortages.

## Nomenclature

**Indices and Sets**

- $t$ — Index of time period, $t = 1, ..., T$.
- $g$ — Index of generator, $g = G1, ..., GN$.
- $\Omega_t$ — Set of time periods.
- $\Omega_g$ — Set of generators.
- $\Omega_{fg}$ — Set of fast-start generators ($\Omega_{fg} \subset \Omega_g$).
- $\Omega_{sg}$ — Set of slow-start generators ($\Omega_{sg} \subset \Omega_g$).

**Parameters**

- $C_g^{NL}$ — No-load cost of generator $g$.
- $C_g^{LP}$ — Linear production cost of generator $g$.
- $C_g^{SU}$ — Startup cost of generator $g$.
- $NL_t$ — Net load for period $t$.
- $P_g^{max}$ — Maximum production capacity of generator $g$.
- $P_g^{min}$ — Minimum production capacity of generator $g$.
- $RR_g$ — Ramp rate of generator $g$.
- $RR_g^{SU}$ — Startup ramp rate of generator $g$.
- $RR_g^{SD}$ — Shutdown ramp rate of generator $g$.
- $P_g^{SUk}$ — Power output of generator $g$ in the $k$th interval of the startup process. Applied only for slow-start units.
- $P_g^{SDk}$ — Power output of generator $g$ in the $k$th interval of the shutdown process. Applied only for slow-start units.
- $SU_g$ — Duration of the startup process of generator $g$. Applied only for slow-start units.
- $SD_g$ — Duration of the shutdown process of generator $g$. Applied only for slow-start units.
- $UFRC_t$ — Requirement for upward ramping capability of system in period $t$.
- $DFRC_t$ — Requirement for downward ramping capability of system in period $t$.
- $\alpha_t$ — Additional ramping capability requirement to cover net load forecasting error.

**Variables**

- $p_{gt}$ — Power output of generator $g$ in period $t$.
- $\overline{p_{gt}}$ — Maximum available power output of generator $g$ in period $t$.
- $ur_{gt}$ — Upward flexible ramping capability of generator $g$ in period $t$.
- $dr_{gt}$ — Downward flexible ramping capability of generator $g$ in period $t$.
- $nur_{gt}$ — Negative contribution of generator $g$ to system upward flexible ramping capability in period $t$.
- $ndr_{gt}$ — Negative contribution of generator $g$ to system downward flexible ramping capability in period $t$.
- $x_{gt}$ — Binary variable that is equal to 1 if generator $g$ is producing above minimum capacity in period $t$ and 0 otherwise.
- $y_{gt}$ — Binary variable that is equal to 1 if generator $g$ starts up in period $t$ 0 otherwise.
- $z_{gt}$ — Binary variable that is equal to 1 if generator $g$ shuts down in period $t$ 0 otherwise.
- $\delta_{gt}^{SD}$ — Auxiliary variable of generator $g$ to calculate negative upward ramping capability in period $t$.
- $\delta_{gt}^{SU}$ — Auxiliary variable of generator $g$ to calculate negative downward ramping capability in period $t$.

This research was supported by Korea Electric Power Corporation. (Grant number: R18XA01).
H. Park, B. Huang, and R. Baldick are with the Department of Electrical and Computer Engineering, The University of Texas at Austin, Austin, TX 78712 USA (e-mail: hyeongon@utexas.edu; binghuang@utexas.edu; baldick@ece.utexas.edu).
H. Park is with the Department of Energy and Electrical Engineering, Korea Polytechnic University, Gyeonggi-do 15073, Korea (e-mail: hyeongon@kpu.ac.kr).

## I. INTRODUCTION

LARGE volumes of renewables, such as solar and wind, are being introduced into electric power systems in order to avoid carbon emissions from fossil fuels and moderate global warming [1]. Increasing penetration of renewable energy is expected to continue owing to recent innovations in renewable generation technologies together with cost competitiveness. However, non-dispatchable characteristics of renewable energy resources is posing new challenges to power system operators: considerable uncertainty and variability in renewable output should be taken into account in system operations [2].

Traditionally, power system operators make unit commitment decisions considering various types of reserves to cover the uncertain and variable nature of net load (*i.e.*, demand minus renewable output). Much attention has been paid to estimate the optimal requirements for reserves to accommodate the increasing amount of renewable energy resources [3]–[5]. This approach, the provision of appropriate amount of deterministic reserves, has advantages in the sense that it could be readily applicable to practical power system operations. Another way to manage uncertainty and variability in net load is the application of a stochastic programming model [6]–[8]. The stochastic programming method attempts to minimize the expected operating costs over the plausible scenarios for net load without any explicit reserve constraints. This approach is economically preferred to the deterministic approach, but it has limitations in computational complexity, defining market settlements, and practical implementation [9]–[11].

Besides reserve products, to further integrate renewables, some independent system operators (ISOs) have introduced ramp products in their electricity markets. The California ISO (CAISO) and the Midcontinent ISO (MISO) have a ramping service known as "fleximamp" and "ramp capability product", respectively, which are designed to improve the operational capability to ramp from one generation level to another over successive dispatch intervals [12], [13]. Flexible ramping capability (FRC) products and reserves share similarities in that both ancillary services set aside a predefined amount of generation capacity. However, the purpose and the expected deployment time of two ancillary services make a clear distinction. FRC products withhold generation capacity, expecting to use the procured capability at one interval later. The main objective of FRC products is to compensate ramping capacity shortage. In contrast, ISOs secure reserves to deal with contingencies that might arise in the same interval.

In recent years, ISOs have already experienced ramping capacity shortage [14], and research on improving FRC in power system operations has become very popular [15]–[22]. Abdul-Rahman *et al.* [15] provided a mathematical formulation of ramp products that can be incorporated with the existing CAISO market rules. Navid *et al.* [16] proposed an FRC model for the MISO market, and derived the cost-effectiveness of the model using various numerical analysis (*e.g.*, simulation on single interval dispatch and time-coupled multi interval dispatch). Optimal requirements for FRC were investigated in [17], where the authors emphasize the reliable operation in high renewable penetrated power systems. The benefits of using wind power generation or electric vehicles as FRC providers were studied in [18], [19]. Wang *et al.* [20] demonstrated the impacts of FRC products on stochastic economic dispatch, and showed market efficiency can be enhanced with FRC products. The authors in [20] extended their research to [21], which focuses on the application of FRC products to unit commitment. In [21], market solutions incorporating FRC products are compared with the benchmark results obtained from stochastic unit commitment. A comprehensive review on the modeling and implementation of FRC products is provided in [22].

However, to the authors' best knowledge, very few publications are available in the literature that address the non-delivery issues of FRC, especially when FRC is considered in unit commitment. In other words, the generation scheduling obtained from the conventional formulations cannot guarantee the availability of FRC even though the solution obtained does satisfy the constraints of FRC requirements. The non-delivery issue, which represents the case where the actual volume of available FRC is less than the calculated FRC volume, might arise when commitment status of generating units is arranged to be changed (*i.e.*, start-up of idle generator, shutdown of online generator). Increasing levels of intermittent renewable generation makes thermal power plants start up and shut down more frequently [23]. It is therefore important to precisely formulate FRC in the unit commitment problem.

This paper aims to reveal the likelihood of ramping capacity shortage in power system operation even if the explicit FRC constraints are considered in the scheduling stage. A new method, herein, has been developed that could manage the non-delivery issue. Our results demonstrate that improperly calculated FRC can be corrected through the proposed method. Moreover, we provide a formulation capable of including startup and shutdown trajectories of slow-start generators in determining FRC. The proposed optimization problem is formulated as a mixed-integer linear programming (MILP) model, which can be solved efficiently using off-the-shelf optimization software.

The remainder of the paper is organized as follows. Section II presents the unit commitment formulation with FRC and addresses the non-delivery issues. The proposed method is provided in Section III. Simulations based on the proposed method are presented in Section IV, and Section V summarizes the results of this work and draws conclusions.

## II. RAMPING PRODUCT IN UNIT COMMITMENT

### A. Conventional Unit Commitment Formulation with Ramping Product

In the CAISO and the MISO markets, two separate types of FRC, termed upward FRC and downward FRC, are co-optimized with energy and other ancillary services. The goal of upward FRC is to alleviate shortage of upward ramping capability, which occurs, for example, when actual output of renewable generation is much smaller than anticipated. On the other hand, downward FRC is secured in preparation for a sudden drop in net load. The requirements for upward FRC and downward FRC, which are calculated just prior to $t = 1$ (*i.e.*,

the beginning point of running real-time look-ahead unit commitment), are given as follows:

$$UFRC_t = \max\left[(NL_{(t+1)} - NL_t) + \alpha_t, 0\right] \quad t=1 \quad (1a)$$

$$DFRC_t = \max\left[(NL_t - NL_{(t+1)}) + \alpha_t, 0\right] \quad t=1 \quad (1b)$$

$$UFRC_t = \max\left[(NL_{(t+1)} - NL_t) + \alpha_t, 0\right] \quad t=2,...,T \quad (1c)$$

$$DFRC_t = \max\left[(NL_t - NL_{(t+1)}) + \alpha_t, 0\right] \quad t=2,...,T \quad (1d)$$

The purpose of upward FRC ($UFRC_t$) is to manage the net load variations between successive intervals (*i.e.*, expected change in net load) and forecast error over the next interval (*i.e.*, unexpected change in net load from what was anticipated). $NL_{(t+1)}$ represents the forecasted net load in time $t+1$, while $NL_t$ represents the current net load. Note that, except for the current net load, all the other net loads are forecasted values, which are anticipated just prior to $t=1$. The volume of FRC required to cover the 15-min ahead net load forecast error is defined as $\alpha_t$.

ISOs schedule the commitment status, generation dispatch, and ancillary services decisions after solving the day-ahead unit commitment problem. The objective of unit commitment is to minimize operating costs while satisfying system constraints (*e.g.*, power balance, reserve provision, and transmission line flow limits) and generating unit constraints (*e.g.*, minimum and maximum output limits, ramp rate limits, minimum on/off time limits). Recently, some ISOs have introduced real-time unit commitment process that determines the commitment status of fast-start generators using 15-minute intervals. The proposed method in this paper can be applied to both day-ahead and real-time unit commitment. However, re-dispatch of units, executed in 5-minute intervals in practical power systems, is not considered for simplicity. Network constraints and other reserve constraints are also neglected to enable a clear interpretation of the obtained results. The mathematical formulation for real-time look-ahead unit commitment with FRC constraints can be modeled as follows [21]:

$$\min \sum_{t \in \Omega_t} \sum_{g \in \Omega_g} \left\{ C_g^{NL} \cdot x_{gt} + C_g^{LP} \cdot p_{gt} + C_g^{SU} \cdot y_{gt} \right\} \quad (2a)$$

Subject to

$$\sum_{g \in \Omega_g} p_{gt} = NL_t, \quad \forall t \quad (2b)$$

$$P_g^{\min} \cdot x_{gt} \leq p_{gt} \leq \overline{p_{gt}} \leq P_g^{\max} \cdot x_{gt}, \quad \forall g,t \quad (2c)$$

$$\overline{p_{gt}} \leq p_{g(t-1)} + RR_g \cdot x_{g(t-1)} + RR_g^{SU} \cdot (x_{gt} - x_{g(t-1)}) + P_g^{\max} \cdot (1 - x_{gt}), \quad \forall g,t \quad (2d)$$

$$\overline{p_{gt}} \leq RR_g^{SD} \cdot (x_{gt} - x_{g(t+1)}) + P_g^{\max} \cdot x_{g(t+1)}, \quad \forall g,t \quad (2e)$$

$$p_{g(t-1)} - p_{gt} \leq RR_g \cdot x_{gt} + RR_g^{SD} \cdot (x_{g(t-1)} - x_{gt}) + P_g^{\max} \cdot (1 - x_{g(t-1)}), \quad \forall g,t \quad (2f)$$

$$P_g^{\min} \cdot (x_{gt} + x_{g(t+1)} - 1) \leq ur_{gt} + p_{gt} \leq \overline{p_{g(t+1)}} + P_g^{\max} \cdot (1 - x_{g(t+1)}), \quad \forall g,t \quad (2g)$$

$$P_g^{\min} \cdot (x_{gt} + x_{g(t+1)} - 1) \leq -dr_{gt} + p_{gt} \leq \overline{p_{g(t+1)}} + P_g^{\max} \cdot (1 - x_{g(t+1)}), \quad \forall g,t \quad (2h)$$

$$-RR_g \cdot x_{g(t+1)} + RR_g^{SD} \cdot (x_{gt} - x_{g(t+1)}) - P_g^{\max} \cdot (1 - x_{gt}) \leq ur_{gt} \leq RR_g \cdot x_{gt}, \quad \forall g,t \quad (2i)$$

$$+RR_g^{SU}(x_{g(t+1)} - x_{gt}) + P_g^{\max} \cdot (1 - x_{g(t+1)}),$$

$$-RR_g \cdot x_{gt} - RR_g^{SU} \cdot (x_{g(t+1)} - x_{gt}) - P_g^{\max} \cdot (1 - x_{g(t+1)}) \leq dr_{gt} \leq RR_g \cdot x_{g(t+1)}, \quad \forall g,t \quad (2j)$$

$$+RR_g^{SD}(x_{gt} - x_{g(t+1)}) + P_g^{\max} \cdot (1 - x_{gt}),$$

$$-P_g^{\max} \cdot x_{gt} + P_g^{\min} \cdot x_{g(t+1)} \leq ur_{gt} \leq P_g^{\max} \cdot x_{g(t+1)}, \quad \forall g,t \quad (2k)$$

$$-P_g^{\max} \cdot x_{g(t+1)} \leq dr_{gt} \leq P_g^{\max} \cdot x_{gt} - P_g^{\min} \cdot x_{g(t+1)}, \quad \forall g,t \quad (2l)$$

$$\sum_{g \in \Omega_g} ur_{gt} \geq UFRC_t, \quad t=1,...,T-1 \quad (2m)$$

$$\sum_{g \in \Omega_g} dr_{gt} \geq DFRC_t, \quad t=1,...,T-1 \quad (2n)$$

$$x_{g(t+1)} - x_{gt} = y_{g(t+1)} - z_{g(t+1)}, \quad \forall g,t \quad (2o)$$

$$x_{gt}, y_{gt}, z_{gt} \in \{0,1\}, \quad \forall g,t \quad (2p)$$

$$p_{gt}, \overline{p_{gt}} \geq 0, \quad \forall g,t \quad (2q)$$

The objective function (2a) is defined to minimize the operating costs which includes generation costs and start-up costs. Equation (2b) represents the power balance that should be maintained, and (2c)–(2f) impose technical limits of each generating unit. The contribution of each unit to FRC is restricted by (2g)–(2l). The constraints relate to upward FRC are formulated as (2g), (2i), and (2k), whereas downward FRC limits are enforced in (2h), (2j), and (2l). The requirements for upward FRC and downward FRC are defined by (2m) and (2n), respectively, of which the threshold values can be determined using (1). The logical constraint for commitment states, startup, and shutdown variables is given in (2o). Binary variables are represented in (2p), and non-negativity constraint applied for power output and available power output of generating units is given by (2q). Note that neither the upward FRC nor downward FRC of generating units is included in (2q).

If trajectories of start-up and shutdown processes of generators are considered in a day-ahead optimization problem, (2b) should be replaced with (3) [24].

$$\sum_{g \in \Omega_g} p_{gt} + \sum_{g \in \Omega_{gs}} \sum_{k=1}^{SU_g} P_g^{SUk} \cdot y_{g,(t-k+SU_g+1)} + \sum_{g \in \Omega_{gs}} \sum_{k=1}^{SD_g} P_g^{SDk} \cdot z_{g,(t-k+1)} = NL_t, \quad \forall t \quad (3)$$

TABLE I. POINTS AT WHICH OPTIMIZATION IS DETERMINED. EMPTY VALUES IN THIS TABLE ARE YET TO BE SOLVED IN PROBLEMS FOR THAT PERIOD.

|  | $t=0$ | $t=1$ | $t=2$ | $t=3$ | $t=4$ | $t=5$ | $t=6$ | $t=7$ |
|---|---|---|---|---|---|---|---|---|
| Prior to $t=1$ | power on/off | on/off |  |  |  |  |  |  |
| In $t=1$ |  | power on/off | power on/off | power on/off | power on/off | power on/off |  |  |
| In $t=2$ |  |  | power on/off | power on/off | power on/off | power on/off | power on/off |  |
| In $t=3$ |  |  | ***power*** ***on/off*** | power *on/off* | power *on/off* | power *on/off* | power *on/off* |  |
| In $t=4$ |  |  |  | power on/off | power on/off | power on/off | power on/off | power on/off |

TABLE II. GENERATOR DATA FOR SIMPLE SYSTEM

| Gen. | $C_g^{LP}$ [$/MWh] | $C_g^{NL}$ [$] | $C_g^{SU}$ [$] | $P_g^{max}$ [MW] | $P_g^{min}$ [MW] | $RR_g$ [MW/15min] | $RR_g^{SU}$ [MW/15min] | $RR_g^{SD}$ [MW/15min] |
|---|---|---|---|---|---|---|---|---|
| G1 | 0 | 0 | 0 | 300 | 300 | 0 | 0 | 0 |
| G2 | 20 | 300 | 300 | 150 | 50 | 40 | 60 | 60 |
| G3 | 40 | 300 | 600 | 200 | 50 | 40 | 60 | 60 |
| G4 | 60 | 300 | 900 | 150 | 50 | 40 | 100 | 100 |

TABLE III. NET LOAD AND RAMP REQUIREMENTS FOR SIMPLE SYSTEM [MW]

|  | $t=1$ | $t=2$ | $t=3$ | $t=4$ | $t=5$ | $t=6$ |
|---|---|---|---|---|---|---|
| Net load ($NL_t$) | **690** | 660 | 640 | 620 |  |  |
| Up ramp ($UFRC_t$) | 0 | 10 | 10 | - |  |  |
| Dn ramp ($DFRC_t$) | 60 | 50 | 50 | - |  |  |
| Net load ($NL_t$) |  | **660** | 640 | 620 | 590 |  |
| Up ramp ($UFRC_t$) |  | 10 | 10 | 0 | - |  |
| Dn ramp ($DFRC_t$) |  | 50 | 50 | 60 | - |  |
| Net load ($NL_t$) |  |  | **665** | 620 | 590 | 570 |
| Up ramp ($UFRC_t$) |  |  | 0 | 0 | 10 | - |
| Dn ramp ($DFRC_t$) |  |  | 75 | 60 | 50 | - |

*B. Non-delivery Issues*

In order to show when and how the non-delivery issue emerges, a simple test system that has four generators is used. Because G1 is the baseload generator whose output is constant and all the other generating units are fast-start units, the trajectory of start-up and shutdown processes are disregarded. In other words, the optimization problem expressed as (2) is solved. The same recursive approach utilized in [21] is adopted to emulate practical power system operations as follows:
1) Prior to $t=1$: Initializes the conditions for the study (*i.e.*, the power output $t=0$ and the commitment decisions for $t=0,1$)
2) In $t=1$: The multi-interval optimization problem is solved to find the on/off states for $t=2,3,4$ and the power output for $t=1,2,3,4$ using the boundary conditions from the first step.
3) In $t=2$: The multi-interval optimization problem is solved to find the on/off states for $t=3,4,5$ and the power output for $t=2,3,4,5$ with the boundary conditions from the second step.
4) In $t=3,4,5,...$: Using the solutions from the previous steps as boundary conditions, the optimization problem is solved. The target time intervals for the commitment decisions and dispatch are rolled forward one interval at a time.

TABLE IV. OPTIMAL SOLUTION EXECUTED IN $t=2$

|  | Product | $t=2$ | $t=3$ | $t=4$ | $t=5$ |
|---|---|---|---|---|---|
| G1 | On/off state ($x_{gt}$) | 1 | 1 | 1 | 1 |
|  | Power ($p_{gt}$) | 300 | 300 | 300 | 300 |
|  | Up ramp ($ur_{gt}$) | 0 | 0 | 0 | - |
|  | Dn ramp ($dr_{gt}$) | 0 | 0 | 0 | - |
| G2 | On/off state ($x_{gt}$) | 1 | 1 | 1 | 1 |
|  | Power ($p_{gt}$) | 150 | 150 | 150 | 150 |
|  | Up ramp ($ur_{gt}$) | 0 | 0 | 0 | - |
|  | Dn ramp ($dr_{gt}$) | 40 | 10~40 | 20~40 | - |
| G3 | On/off state ($x_{gt}$) | 1 | 1 | 1 | 1 |
|  | Power ($p_{gt}$) | 160 | 190 | 170 | 140 |
|  | Up ramp ($ur_{gt}$) | 10~40 | 10 | 0~30 | - |
|  | Dn ramp ($dr_{gt}$) | 0~40 | 10~40 | 20~40 | - |
| G4 | On/off state ($x_{gt}$) | 1 | 0 | 0 | 0 |
|  | Power ($p_{gt}$) | 50 | 0 | 0 | 0 |
|  | Up ramp ($ur_{gt}$) | 0 | 0 | 0 | - |
|  | Dn ramp ($dr_{gt}$) | 50 | 0 | 0 | - |
|  | Operating cost [k$] | 13.3 | 11.2 | 10.4 | 9.2 |

Table I shows the determining points at which optimization is calculated. The shaded parts represent the variables that are already determined in that period, which are used as boundary constraints. For example, variables that are boldly italicized represent the boundary conditions - determined in $t=1,2$ - for the optimization problem executed in $t=3$. The italic variables are the variables to be optimized in $t=3$.

The specifications for units and input data for the system are tabulated in Table II and Table III, respectively. It is assumed that net load can vary in the range of ±30MW, so FRC requirements for covering this uncertainty ($\alpha_t$) is set to 30MW. Thus, for example, upward FRC requirements in $t=2$ can be calculated as (640MW-660MW)+30MW = 10MW according to (1). Likewise, it can be computed that 50MW downward FRC is needed in $t=2$. Note that the first values in net load rows (italic and bold numbers) represent realized net load for each period, whereas the other net load values are forecasted values.

Table IV shows the optimal solution executed in $t=2$ based on the boundary condition obtained from the solution in $t=1$. The optimal values of procured volume of FRC are represented as continuous ranges if the constraints for FRC requirements are not binding (*e.g.*, $ur_{gt}$ and $dr_{gt}$ of G3 in $t=2$). Any combinations of the values meeting the requirements yield the same objective value.

The generator G3 can ramp up as much as 40 MW from $t=2$ to $t=3$, therefore, the maximum ramp up capability of the system in $t=2$ is 40 MW. However, it must be noted that the system cannot actually ramp up 40MW until $t=3$. The actual "deliverable" upward FRC procured in $t=2$ is -10MW. The reason for this miscalculation can be explained as follows.

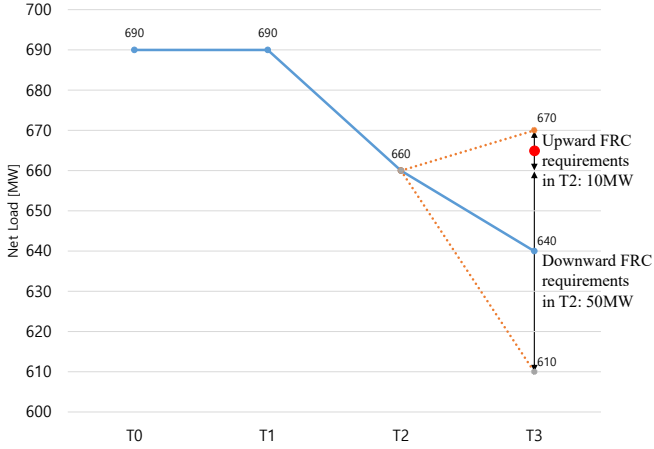

Fig. 1: Net load forecast marked in blue solid line and FRC requirements in $t=2$. Red dot represents the realized net load in $t=3$.

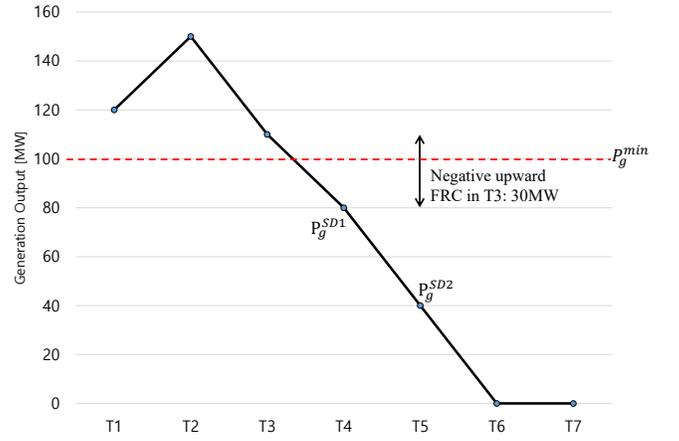

Fig. 2: Shutdown trajectory of slow-start generating unit g. The unit is online from $t=1$ to $t=5$ and becomes offline from T6. $x=1$ from $t=1$ to $t=3$ and 0 otherwise.

TABLE V. OPTIMAL SOLUTION EXECUTED IN $t=3$ AND OPERATING COSTS

|  | Product | $t=3$ | $t=4$ | $t=5$ | $t=6$ |
|---|---|---|---|---|---|
| G1 | Power ($p_{gt}$) | 300 | 300 | 300 | 300 |
| G2 | Power ($p_{gt}$) | 150 | 150 | 150 | 150 |
| G3 | Power ($p_{gt}$) | 200 | 170 | 140 | 120 |
| G4 | Power ($p_{gt}$) | 0 | 0 | 0 | 0 |
| Load shedding [MW] | | 15 | 0 | 0 | 0 |
| Net load [MW] | | 665 (realized) | 620 | 590 | 570 |
| Operating cost [k$] | | 146.6 | 10.4 | 9.2 | 9.6 |

In $t=2$, when making the commitment decisions for $t=3$, the generating unit G4 is determined to shut down in $t=3$, which results in zero power output of G4 in $t=3$. The problem is that the drop of 50MW in output from G4 (the corresponding volume that G4 produces in $t=2$) curtails system upward FRC because other generators should ramp up as much as 50MW to make up for the output decrease of G4. The upward FRC contribution of G4 at $t=2$ should be calculated as -50MW. However, the obtained solution by solving (2) indicates that there is zero upward FRC contribution of G4 in $t=2$. It is interesting to note that although non-negativity constraint (2r) exclude the FRC contributions of generating units, it is not enough to account for the "negative" effect of G4 on upward FRC in the system. In order to clarify the "negative" effects of the generating units scheduled to be turned on or off, we define new variables $nur_{gt}$ ($ndr_{gt}$) to represent negative contribution to system upward (downward) FRC from the generator $g$ in period $t$.

Table V shows the optimization results executed in $t=3$ when realized net load is 655MW which lies in the range of 610MW to 670MW. Even though the realized uncertainty is +25MW, which is within the minimum and maximum bounds (Fig. 1), it can be seen from the table that involuntary load shedding occurs in $t=3$ due to the non-delivery issue. The operating costs for $t=3$ can be computed as $146.6 k if the value of lost load (VOLL) is $9,000/MWh.

## III. PROPOSED APPROACH

### A. Additional Constraints for Start-up Generators and Shutdown Generators

We derive new constraints to take into account the reduced FRC because of the generating units that are planned to be turned on or off. We adopt the sign convention that a positive value of $nur_{gt}$ represents a negative contribution to upward FRC. Likewise, a new variable $ndr_{gt}$ is defined to represent a negative effect to downward FRC. The following formulations (4) can deal with the non-delivery issue.

$$p_{gt} = nur_{gt} + \delta_{gt}^{SD}, \quad \forall g \in \Omega_{fg}, t \quad (4a)$$

$$P_g^{\min} \cdot z_{g(t+1)} \leq nur_{gt} \leq P_g^{\max} \cdot z_{g(t+1)}, \quad \forall g \in \Omega_{fg}, t \quad (4b)$$

$$0 \leq \delta_{gt}^{SD} \leq P_g^{\max} \cdot (1 - z_{g(t+1)}), \quad \forall g \in \Omega_{fg}, t \quad (4c)$$

$$p_{gt} = ndr_{gt} + \delta_{gt}^{SU}, \quad \forall g \in \Omega_{fg}, t \quad (4d)$$

$$P_g^{\min} \cdot y_{g(t+1)} \leq ndr_{gt} \leq P_g^{\max} \cdot y_{g(t+1)}, \quad \forall g \in \Omega_{fg}, t \quad (4e)$$

$$0 \leq \delta_{gt}^{SU} \leq P_g^{\max} \cdot (1 - y_{g(t+1)}), \quad \forall g \in \Omega_{fg}, t \quad (4f)$$

Constraints (4a)–(4c) are applied for negative upward FRC, and constraints (4d)–(4f) are for negative downward FRC. If the generator $g$ that produces energy above the minimum output level in $t$ is turned off in $t+1$, the subsidiary variable $\delta_{gt}^{SD}$ becomes zero while the binary variable $z_{g(t+1)}$ has a value of one. In this case, constraints (4a)–(4c) enforce that the volume of the negative upward FRC should be equal to the generation output in $t$ ($p_{gt}$). In all the other cases (i.e., a generating unit continues to generate power above the minimum power limit until $t+1$ or a unit is offline in $t$), $z_{g(t+1)}$ has a value of zero, which results in a zero negative effect on downward FRC ($nur_{gt}=0$). Similarly, if the generating unit is scheduled to start up, negative downward FRC has a nonzero value, in other words, it curtails the downward FRC of the system.

Generating units that complete startup and shutdown process within one interval can be modeled as (4). However, if the procedure of startup and shutdown of a generator takes longer than one interval, which we referred to as a slow-start generator, different formulations should be derived to reflect startup and shutdown trajectories. Fig. 2 illustrates an example trajectory of a slow-start generator whose duration of the shutdown process is two intervals. The red dashed line in Fig. 2 represents the minimum capacity of the generating unit. If the generating unit is scheduled to be turned off in $t = 4$, this unit is in the process of shutdown for two intervals which leads the system operator to consider the reduced ramping capability for three intervals.

A negative upward FRC of the generating unit in $t = 3$, which is the last period that the unit generates above the minimum power limit, depends on the generation output of the unit itself ($p_{g,T3}$). The negative upward FRC is therefore a continuous variable, which can be computed as $p_{g,T3} - P_g^{SD1}$. Here, $P_g^{SD1}$ represents the output of the first segment of the shutdown trajectory. By employing the same technique used in formulating (4), that is by introducing auxiliary variables, negative upward FRC can be formulated as (5a)–(5c). On the other hand, negative upward FRC in $t = 4$ and $t = 5$ have the fixed values, $P_g^{SD1} - P_g^{SD2}$ and $P_g^{SD2}$, respectively, because the power output in the shutdown process are predefined constant values. The formulation of negative upward FRC applied for the generating unit whose duration of the shutdown process is $SD_g$ can be expressed as follows:

$$p_{gt} - P_g^{SD1} \cdot x_{gt} = nur_{gt} + \delta_{gt}^{SD} + P_g^{SD1}, \quad \forall g \in \Omega_{sg}, t \quad (5a)$$

$$(P_g^{\min} - P_g^{SD1}) \cdot z_{g(t+1)} \leq nur_{gt}$$
$$\leq (P_g^{\max} - P_g^{SD1}) \cdot z_{g(t+1)}, \quad \forall g \in \Omega_{sg}, t \quad (5b)$$

$$0 \leq \delta_{gt}^{SD} + P_g^{SD1} \leq P_g^{\max} \cdot (1 - z_{g(t+1)}), \quad \forall g \in \Omega_{sg}, t \quad (5c)$$

$$nur_{gt}^{SD} = (P_g^{SD1} - P_g^{SD2}) \cdot z_{gt}$$
$$+ (P_g^{SD2} - P_g^{SD3}) \cdot z_{g(t-1)} + \ldots \quad \forall g \in \Omega_{sg}, t \quad (5d)$$
$$+ P_g^{SD(SD_g)} \cdot z_{g(t-SD_g+1)},$$

Similarly, the negative downward FRC can be formulated as follows:

$$p_{g(t+1)} - P_g^{SU(SU_g)} \cdot x_{g(t+1)}$$
$$= ndr_{gt} + \delta_{gt}^{SU} + P_g^{SU(SU_g)}, \quad \forall g \in \Omega_{sg}, t \quad (5e)$$

$$(P_g^{\min} - P_g^{SU(SU_g)}) \cdot y_{g(t+1)} \leq ndr_{gt}$$
$$\leq (P_g^{\max} - P_g^{SU(SU_g)}) \cdot y_{g(t+1)}, \quad \forall g \in \Omega_{sg}, t \quad (5f)$$

$$0 \leq \delta_{gt}^{SU} + P_g^{SU(SU_g)} \leq P_g^{\max} \cdot (1 - y_{g(t+1)}), \quad \forall g \in \Omega_{sg}, t \quad (5g)$$

$$ndr_{gt}^{SU} = (P_g^{SU(SU_g)} - P_g^{SU(SU_g-1)}) \cdot y_{g(t+2)}$$
$$+ (P_g^{SU(SU_g-1)} - P_g^{SU(SU_g-2)}) \cdot y_{g(t+3)} \quad \forall g \in \Omega_{sg}, t \quad (5h)$$
$$+ \ldots + P_g^{SU1} \cdot y_{g(t+SU_g+1)},$$

### B. Reformulation of Existing Constraint

In order to apply the proposed modeling of negative FRC, the upward and downward ramping requirement constraints, should be reformulated as (6a) and (6b), respectively.

$$\sum_{g \in \Omega_g} ur_{gt} - \sum_{g \in \Omega_g} nur_{gt} - \sum_{g \in \Omega_{sg}} nur_{gt}^{SD} \geq UFRC_t, \quad \forall t \quad (6a)$$

$$\sum_{g \in \Omega_g} dr_{gt} - \sum_{g \in \Omega_g} ndr_{gt} - \sum_{g \in \Omega_{sg}} ndr_{gt}^{SU} \geq DFRC_t, \quad \forall t \quad (6b)$$

It should be noted that the optimal FRC requirements for the system is not modified, nor is it calculated in a different way. (Compare (6a) and (6b) with (2h) and (2i). In this paper, we derive new formulations that enable the "deliverable" ramping capacity to be obtained if the optimal ramping requirements are given. The optimal requirements can be improved through various methods [17], and the proposed method can be suitably applied once the requirements are known in advance.

TABLE VI. OPTIMAL SOLUTION EXECUTED IN $t = 2$ WITH ADDITIONAL CONSTRAINTS (PR.: PROPOSED METHOD, CO.: CONVENTIONAL METHOD)

|  | Product | T2 | | T3 | | T4 | | T5 | |
| --- | --- | --- | --- | --- | --- | --- | --- | --- | --- |
|  |  | Pr. | Co. | Pr. | Co. | Pr. | Co. | Pr. | Co. |
| G1 | On/off state | 1 | 1 | 1 | 1 | 1 | 1 | 1 | 1 |
|  | Power | 300 | 300 | 300 | 300 | 300 | 300 | 300 | 300 |
|  | Up ramp | 0 | 0 | 0 | 0 | 0 | 0 | - | - |
|  | Dn ramp | 0 | 0 | 0 | 0 | 0 | 0 | - | - |
| G2 | On/off state | 1 | 1 | 1 | 1 | 1 | 1 | 1 | 1 |
|  | Power | 150 | 150 | 130 | 150 | 150 | 150 | 150 | 150 |
|  | Up ramp | 0 | 0 | 20 | 0 | 0 | 0 | - | - |
|  | Dn ramp | 40 | 40 | 40 | 40 | 20 | 40 | - | - |
| G3 | On/off state | 1 | 1 | 1 | 1 | 1 | 1 | 1 | 1 |
|  | Power | 160 | 160 | 160 | 190 | 170 | 170 | 140 | 140 |
|  | Up ramp | 40 | 40 | 40 | 10 | 30 | 30 | - | - |
|  | Dn ramp | 40 | 40 | 40 | 40 | 40 | 40 | - | - |
| G4 | On/off state | 1 | 1 | **1** | **0** | 0 | 0 | 0 | 0 |
|  | Power | 50 | 50 | 50 | 0 | 0 | 0 | 0 | 0 |
|  | Up ramp | 40 | 0 | 0 | 0 | 0 | 0 | - | - |
|  | Dn ramp | 0 | 50 | 50 | 0 | 0 | 0 | - | - |
| Operating cost [k$] | | 13.3 | | 12.9 | | 10.4 | | 9.2 | |

TABLE VII. OPTIMAL SOLUTION EXECUTED IN $t = 3$ AND OPERATING COSTS BASED ON PROPOSED METHOD

| | Product | T3 | T4 | T5 | T6 |
| --- | --- | --- | --- | --- | --- |
| G1 | Power ($p_{gt}$) | 300 | 300 | 300 | 300 |
| G2 | Power ($p_{gt}$) | 150 | 140 | 150 | 150 |
| G3 | Power ($p_{gt}$) | 165 | 130 | 140 | 120 |
| G4 | Power ($p_{gt}$) | 50 | 50 | 0 | 0 |
| Load shedding | | 0 | 0 | 0 | 0 |
| Net load | | 655 (realized) | 620 | 590 | 570 |
| Operating cost [k$] | | 13.5 | 11.9 | 9.2 | 8.4 |

## IV. TEST RESULTS

The proposed formulation is tested for the simple system and the modified IEEE 118-bus system. The same simple problem introduced in Section II-B is analyzed to show that the non-delivery issue can be avoided with the proposed approach. In order to test the scalability of the proposed method, the day-ahead unit commitment problem is solved based on the IEEE 118-bus system which includes wind power generators. All simulations in this paper were carried out on a personal computer with a 3.60-GHz Intel Core i5 8600K CPU, 16-GB RAM, and 64-bit operating system. The optimization solver GUROBI under GAMS was used to solve the problem, and the relative optimality tolerance was set to 0.1%.

### A. Toy Example

The performance of the new method is evaluated on the system introduced in Section II-B. The optimization problem that comprises the objective function (2a) and the constraints (2b)-(2r) and (4) is solved based on the same input data. Table VI compares the optimal solution with the proposed method and with the conventional method. Among the optimal solutions, the case where the procured volume of FRC can be maximized is listed in the table. The changes in the unit commitment results executed in $t=2$ are listed as underlined data in Table VI. Due to the new constraints to deal with the non-delivery issue, G4 stayed online until $t=3$ in the proposed case while the commitment schedules of the other units were the same as those of the conventional case.

When the realized net load in T3 is 655MW, which is the same load level as the example in Section II-B, the optimal solution executed in $t=3$ can be obtained as is summarized in Table VII. As can be seen in Table VII, no involuntary load shedding occurs in $t=3$ due to the sufficient and reliable ramping capability procured in $t=2$. It should be remembered that load shedding may take place if the unit commitment is solved with the conventional approach (Compare Table V).

### B. Modified IEEE 118-Bus System

The proposed method was applied to a more realistic problem with the modified IEEE 118-bus test system, which has 54 slow-start generating units. The total installed wind capacity was assumed to be 50% of peak load. The forecasted values of the demand and wind power generation for each hour is listed in Table VIII. The specific data of units were taken from [24].

The day-ahead FRC requirements are determined based on estimates of the real-time ramping needs, which will not be known until the operating day. In this study, the day-ahead FRC requirements are computed similar to the real-time case, as expressed in (1). The hourly variabilities of the net load are considered, and hour ahead forecast error distributions are used instead of using 15-min ahead forecast error. It is assumed that both the demand forecasting error and the wind power generation forecasting error follow a normal distribution with zero mean [25]. The standard deviation for the demand forecasting error was set to 1% of the forecasted demand, and the standard deviation for the wind power generation error was set to 4% of the installed wind capacity [26].

TABLE VIII. FORECASTED VALUES OF THE DEMAND AND WIND POWER GENERATION

| Time | 1 | 2 | 3 | 4 | 5 | 6 |
|---|---|---|---|---|---|---|
| Demand [MW] | 4920 | 3960 | 3480 | 2400 | 3000 | 3600 |
| Wind [MW] | 300 | 292.5 | 307.5 | 315 | 285 | 277.5 |
| Time | 7 | 8 | 9 | 10 | 11 | 12 |
| Demand [MW] | 4200 | 4680 | 4920 | 5280 | 5340 | 5040 |
| Wind [MW] | 285 | 292.5 | 262.5 | 247.5 | 255 | 292.5 |
| Time | 13 | 14 | 15 | 16 | 17 | 18 |
| Demand [MW] | 4800 | 4560 | 5280 | 5400 | 5100 | 5340 |
| Wind [MW] | 307.5 | 322.5 | 307.5 | 285 | 262.5 | 240 |
| Time | 19 | 20 | 21 | 22 | 23 | 24 |
| Demand [MW] | 5640 | 5880 | 6000 | 5400 | 5220 | 4920 |
| Wind [MW] | 225 | 217.5 | 240 | 255 | 262.5 | 247.5 |

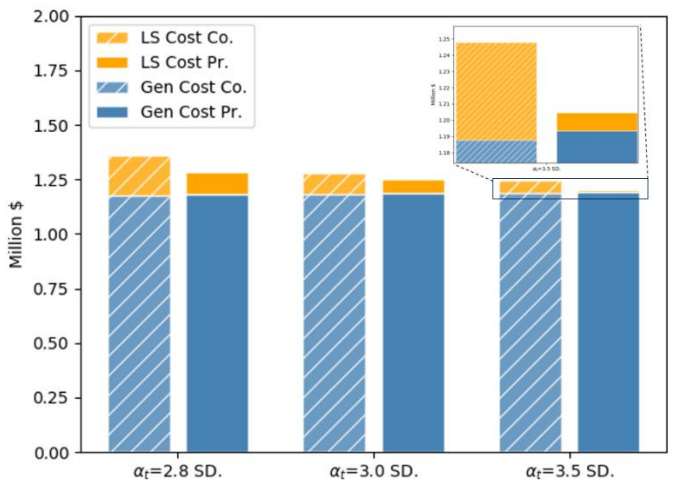

Fig. 3: Evaluation results of the proposed and conventional unit commitment solutions with 2500 validation scenarios for the IEEE 118-bus system. (VOLL=$9,000/MWh, Pr.: proposed method, Co.: conventional method)

In order to compare the conventional approach and the proposed approach, we generated 2500 scenarios using the Monte Carlo simulation. The performance of each approach is analyzed using the following procedure, as was done in [27]:

1) Solve the hourly unit commitment problem with 54 units for the 24-hour periods under the central forecast net load. The conventional approach uses (2a), (2c)-(2q), (3), while (2a), (2c)–(2q), (3)-(6) are employed for the proposed approach. Minimum on/off time limits, which are not shown herein, are also considered when solving the problem.
2) Select one of the generated scenarios.
3) Evaluate the performance of both the conventional unit commitment solution and the proposed unit commitment solution that are obtained in the first step with the selected scenario. The online generating units are re-dispatched to cover the realized uncertainty, and the operating cost and load-shedding cost are calculated for each of the approaches. If load shedding is needed to satisfy the power balance constraint, this is penalized by a cost in the objective function.
4) Go to the third step above using another scenario and repeat the process until the last scenario (the 2500[th] scenario). Compute and save the average values of the generation cost and the load-shedding cost.
5) Compare the results. The expected operating cost is defined as the sum of the average generation cost and the average load-shedding cost.

Fig. 3 summarizes the expected operating costs of the proposed and the conventional unit commitment solutions. The volume of additional FRC requirements to cover forecasting error ($\alpha_t$) was set to 2.8, 3.0, and 3.5 times the standard deviation of the forecasting error of net load. If the procured volume of FRC is not enough, involuntary load-shedding events occur even in the proposed method. However, it should be noted that the proposed method can ensure a reliable operation for any deviation which lies within the target bounds. With the conventional unit commitment solution, we have found some cases, similar to the toy-example, where the SO needs to resort to load shedding even though the deviation of the forecasted value is within the minimum and maximum error bounds.

The computing times required to solve the problem are listed in Table IX. The additional variables and constraints increase the computational burden, but the solution times in all cases are still less than 20 min, which is a short enough time in the day-ahead scheduling [27].

TABLE IX. COMPARISON OF COMPUTATIONAL TIME [S]

| $\alpha$ | 2.8 times std. dev. | 3.0 times std. dev. | 3.5 times std. dev. |
|---|---|---|---|
| Proposed method | 947 | 1107 | 898 |
| Conventional method | 745 | 702 | 831 |

## V. CONCLUSION

In this study, we proposed a new formulation that can ensure the deliverability of the procured FRC. A method to consider the trajectory of the startup and shutdown process of the generator in the determination of FRC was also developed. Newly derived formulations have been expressed in the MILP form to be easily solved.

Simulation results have shown that a generation scheduling, based on conventional formulation, might yield unreliable operations. We have found cases where the procured volume of FRC is insufficient, even if the deviation of the forecasted values is within the anticipated bounds, which leads to load shedding. With the proposed model, an SO can guarantee the reliable operations if the optimal requirements for FRC are properly predefined.

The specifications of ramping products and the requirements for the volume of FRC may differ from case to case. However, we believe that the proposed model can generally be applied to other scenarios that focus on the design of the ramping products.